\documentclass[letterpaper, dvipsnames, 12pt]{article}
\usepackage[margin=1in]{geometry}
\usepackage[utf8]{inputenc}
\usepackage[T1]{fontenc}
\usepackage{amsmath,amsthm,amssymb}

\usepackage{graphicx}
\usepackage{tikz}
\usepackage{hyperref}
\usepackage{color}
\usepackage{enumerate}
\usepackage{multicol}
\usepackage{xcolor}
\usepackage{amsrefs}
\newcommand{\SL}{\sum\limits_}

\newcommand{\R}{\mathbb{R}}

\newcommand{\PR}{\mathbb{P}}
\newcommand{\E}{\mathbb{E}}

\usepackage[autostyle,english=american]{csquotes}
\MakeOuterQuote{"}

\newtheorem{ccounter}{ccounter}[section]
\newtheorem{thm}[ccounter]{Theorem}
\newtheorem{lem}[ccounter]{Lemma}

\newtheorem{defn}[ccounter]{Definition}
\newtheorem{prop}[ccounter]{Proposition}

\title{Polynomial Mixing of the critical Glauber Dynamics for the Ising Model}
\author{Kyprianos-Iason Prodromidis\thanks{Department of Mathematics, Princeton University, Email: kp2702@princeton.edu} \and Allan Sly\thanks{Department of Mathematics, Princeton University, Email: allansly@princeton.edu}}
\date{}
\begin{document}
\maketitle
\begin{abstract}
In this note, we prove that on any graph of maximal degree $d$ the mixing time of the Glauber Dynamics for the Ising Model at $\beta_c=\tanh^{-1}(\frac1{d-1})$, the uniqueness threshold on the infinite $d$-regular tree, is at most polynomial in $n$.  The proof follows by a simple combination of new log-Sobolev bounds of Bauerschmidt and Dagallier, together with the tree of self avoiding walks construction of Weitz.
 While preparing this note we became aware that Chen, Chen, Yin and Zhang recently posted another proof of this result. We believe the simplicity of our argument is of independent interest.
\end{abstract}
\section{Introduction}
Weitz~\cite{Weitz} showed that on any $d$-regular graph the Ising model has spatial decay of correlation for any inverse temperature $|\beta|<\beta_c$ where
\[
(d-1)\tanh(\beta_c)=1,
\]
is the uniqueness threshold for the infinite $d$-regular tree.   In particular, the tree is extremal in terms of correlation decay.  A common folklore principle in mixing times of Glauber dynamics is that spatial decay of correlation should imply fast mixing of the Glauber dynamics.  Here we consider the critical case and, on an $n$-vertex graph $G$ with maximal degree $d$, we show that the Glauber Dynamics at  $\beta_c$ mixes in polynomial time.
\begin{thm}\label{main thm}
Let $d\ge3$ be a positive integer. There exists a universal constant $C>0$, independent of $d$, such that for $n$ large enough, if $G$ is a $d$-regular graph on $n$ vertices, then the Glauber dynamics for the critical Ising model on $G$ has
$$t_{\beta_c,\text{mix}}\le n^C.$$
\end{thm}

In 2007, Gerschenfeld and Montanari calculated the free energy of the Ising model on random $d$-regular graphs~\cite{G-M} which implies that in the low-temperature regime $\beta>\beta_c$, the mixing time is $\exp(\Omega(n))$. By contrast in 2009,  Mossel and the second author~\cite{M-S} studied the case in which $\beta<\beta_c$, and proved that the mixing time is of $O_\beta(n\log(n))$. The proof combined the censoring inequality of Peres and Winkler~\cite{P-W} with estimates from the tree of Self-Avoiding Walks construction introduced by Weitz~\cite{Weitz} to iteratively control the magnetization.

With the advent of Stochastic Localization techniques, new methods of studying the mixing time have been developed, including \cite{C-E} and \cite{B-D} established explicit upper bounds for the log-Sobolev constants of Glauber Dynamics. Most relevantly~\cite{B-D} showed a polynomial upper bound on the mixing time at the critical temperature for lattice boxes in $\mathbb{Z}^d$ when $d\geq 5$ using the fact that the 0-field susceptibility is finite.

Shortly before posting this work, we learned that Chen, Chen, Yin and Zhang~\cite{C-C-Y-Z}, studied the same Markov chain as Theorem~\ref{main thm} establishing an upper bound of $O(n^{2+\frac{2}{d-2}})$ for the mixing time, which gives a better exponent. Their paper contains numerous additional results including the analogous result for the hardcore model as well as polynomial lower bounds for certain random graphs. Their approach uses a Theorem from~\cite{C-E} together with new spectral independence estimates. Our approach, using the main theorem of \cite{B-D}, requires only basic estimates of correlations and we believe the straightforward nature of our analysis is noteworthy in its own right.
\section{Proof of the Main Theorem}
\subsection{Basic Definitions}
In this subsection, we give the definitions and reduce the problem to the bounding of the quantity $\chi_s$, defined in Theorem \ref{Ineq-Log-Sobolev-Constant}.
\begin{defn}
\begin{itemize}
\item 
Let $G$ be a graph on a vertex set $V$. The Ising model on this graph with inverse temperature $\beta$ and external field $h:V\to[-\infty,\infty]$ is the probability measure $\mu$ on $\{-1,1\}^V$ satisfying
$$\mu_{\beta,h}(\sigma)=\dfrac{1}{Z_{G,h}}\exp\left(\beta\SL{u\sim v}\sigma_u\sigma_v+\SL{v\in V}h(v)\sigma_v\right).$$
By $\mu_\beta$ we will just denote the Ising model with $h=0$.
\item
For a configuration $\sigma\in\{-1,1\}^V$ and for some $i\in[n]$, let $\sigma^i$ be the configuration $\sigma$ in which we have flipped coordinate $i$. For two configurations $\sigma,\tau\in\{-1,1\}^V$, we write $\sigma\sim\tau$ if $\tau=\sigma^i$ for some $i\in[n]$. The Glauber Dynamics on $G$ is a Markov chain with transition probabilities, for $\sigma\neq\tau$, equal to
$$P_\beta(\sigma,\tau)=\begin{cases}
\dfrac{\mu_\beta(\tau)}{\mu_\beta(\sigma)+\mu_\beta(\tau)},\ \ \text{if}\ \ \sigma \sim \tau
\\ \\
0,\ \ \text{otherwise}.
\end{cases}$$
\item
On the set of probability measures on $\{-1,1\}^n$, we define the \textbf{Total Variation Distance} between $\mu,\nu$ to be
$$\text{d}_{\text{TV}}(\mu,\nu):=\sup\limits_{A}|\mu(A)-\nu(A)|=\dfrac{1}{2}\SL{\sigma}|\mu(\sigma)-\nu(\sigma)|.$$
\item 
The mixing time of this Markov Chain is defined to be
$$t_{\beta,\text{mix}}(\varepsilon)=\inf\{t\ge0:\text{d}_{\text{TV}}(P_\beta^t\delta_{\sigma_0},\mu_\beta)\le\varepsilon,\ \forall\sigma_0\in\{-1,1\}^V\},$$
where $\delta_{\sigma_0}$ is the probability measure that gives mass 1 to $\sigma_0$. We set
$t_{\beta,\text{mix}}:=t_{\beta,\text{mix}}\left(1/4\right).$
\end{itemize}
\end{defn}
\begin{defn}
Let $f:\{-1,1\}^V\to[0,\infty)$ be a function on the possible configurations. We define the \textbf{entropy} of $f$ to be
$$\text{Ent}_{\mu_\beta}(f):=\E_{\mu_\beta}(f\log(f))-\E_{\mu_\beta}(f)\log\E_{\mu_\beta}(f),$$
where we set $0\log0=0$.
Also, for a function $f:\{-1,1\}^V\to\R$, we define the Dirichlet forms $\mathcal{E}_1$ and $\mathcal{E}_2$ as follows:
$$\mathcal{E}_1(f,f):=\SL{i=1}^n\E_\beta\left[(f(\sigma)-f(\sigma^i))^2\right],\ \ \text{and}\ \ \mathcal{E}_2(f,f):=\dfrac{1}{2}\SL{\sigma,\tau}\mu_\beta(\sigma)P_\beta(\sigma,\tau)(f(\sigma)-f(\tau))^2.$$
Also, let $\gamma_\beta$ and $\alpha_\beta$ be the respective log-Sobolev constants, i.e.
$$\gamma_\beta=\inf\left\{\dfrac{\mathcal{E}_1(f,f)}{\text{Ent}_{\mu_\beta}(f^2)}:\text{Ent}_{\mu_\beta}(f^2)\neq0\right\}\ \ \text{and}\ \ \alpha_\beta=\inf\left\{\dfrac{\mathcal{E}_2(f,f)}{\text{Ent}_{\mu_\beta}(f^2)}:\text{Ent}_{\mu_\beta}(f^2)\neq0\right\}.$$
\end{defn}
In the proof of Theorem \ref{main thm}, we will use the following (which is a special case of the main Theorem proven in \cite{B-D}).
\begin{thm}\label{Ineq-Log-Sobolev-Constant}
For the Ising model given by
$$\PR(\sigma)\varpropto\exp\left(\dfrac{\beta}{2}\langle\sigma,A\sigma\rangle\right),$$
the log-Sobolev constant $\gamma_\beta$ satisfies the inequality
$$\gamma_\beta^{-1}\le\dfrac{1}{2}+\beta\cdot\lVert A\rVert\cdot\exp\left(2\lVert A\rVert\int_0^\beta\chi_s\text{d}s\right),$$
where $\lVert A\rVert$ is the spectral radius of $A$, $\chi_s:=\sup\limits_{x\in V(G)}\SL{y\in V(G)}\E_s(\sigma_x\sigma_y)$ and $\E_s$ denotes the expected value with inverse temperature equal to $s$.
\end{thm}
We also state a Theorem which relates the log-Sobolev constant $\alpha_\beta$ with the mixing time of $P_\beta$, and in fact, holds for general Markov chains (for a proof, we refer to \cite{D-SC}):
\begin{thm}\label{log-Sobolev and mixing}
For the Ising model on $G$, the following inequality holds:
$$t_{\text{mix}}(\varepsilon)\le\dfrac{1}{4\alpha_\beta}\left(\log\log\left(\dfrac{1}{\mu_{\beta,\min}}\right)+\log\left(\dfrac{1}{2\varepsilon^2}\right)\right).$$
\end{thm}
We now explain how the two log-Sobolev constants compare to each other. Because $G$ is a $d$-regular graph, for $\sigma\sim\tau$ we always have
$$e^{-2d\beta}\le\dfrac{\mu_\beta(\sigma)}{\mu_\beta(\tau)}\le e^{2d\beta},$$
therefore
$$P_\beta(\sigma,\tau)=\dfrac{1}{n}\cdot\dfrac{\mu_\beta(\tau)}{\mu_\beta(\sigma)+\mu_\beta(\tau)}\ge\dfrac{1}{n}\cdot\dfrac{1}{1+e^{2d\beta}}.$$
This implies that
$$\mathcal{E}_2(f,f)\ge\dfrac{1}{2n(1+e^{2d\beta })}\cdot\SL{\sigma,\tau}\mu_\beta(\sigma)(f(\sigma)-f(\tau))^2=\dfrac{1}{2n(1+e^{2d\beta })}\cdot\mathcal{E}_1(f,f),$$
so
\begin{align}\label{comparison-log-Sobolev-constants}
\alpha_\beta\ge\dfrac{1}{2n(1+e^{2d\beta})}\cdot\gamma_\beta\ \ \Rightarrow\ \ \alpha_\beta^{-1}\le 200\cdot n\cdot\gamma_\beta^{-1}.
\end{align}
\begin{lem}\label{main-lem}
For the quantity $\chi_s$, the inequality
$$\chi_s\le\min\left(\dfrac{100d}{(d-1)^2(\beta_c-s)},n\right)$$
is true.
\end{lem}
We explain why this Lemma finishes the proof of Theorem \ref{main thm}.
\begin{proof}[Proof of Theorem \ref{main thm} given Lemma \ref{main-lem}]
Note that because of Theorem \ref{log-Sobolev and mixing}, it suffices to prove that
$$\alpha_{\beta_c}^{-1}\le 100\cdot n^{1+\frac{200d^2}{(d-1)^2}}.$$
We use Lemma \ref{main-lem}:
$$\int_0^{\beta_c}\chi_s\ \text{d}s\le\dfrac{100d}{(d-1)^2}\int_0^{\beta_c-\frac{1}{n}}\dfrac{1}{\beta_c-s}\ \text{d}s+\int_{\beta_c-\frac{1}{n}}^{\beta_c}n\ \text{d}s=\dfrac{100d}{(d-1)^2}\cdot(\log(n)-\log(\beta_c))+1.$$
Using the inequality of Theorem \ref{Ineq-Log-Sobolev-Constant}, it is clear that
$$\gamma_{\beta_c}^{-1}\le Cn^{\frac{200d^2}{(d-1)^2}},$$
and because of relation (\ref{comparison-log-Sobolev-constants}), the proof of Theorem \ref{main thm} is complete.
\end{proof}
\subsection{Proof of Lemma \ref{main-lem}}
In this section, we analyze the sum of correlations that appears in Theorem \ref{Ineq-Log-Sobolev-Constant}, and combine it with the tree of Self-Avoiding Walks construction, to conclude the proof of Lemma \ref{main-lem}, and thus the proof of Theorem \ref{main thm}.

Since the inequality $\chi_s\le n$ is obvious, we will only prove that $\chi_s\le\dfrac{100d}{(d-1)^2(\beta_c-s)}.$ The main tool towards proving that is the following Proposition:
\begin{prop}\label{main-prop}
Let $\mathcal{T}$ be a tree and let $(\sigma_v)_{v\in V(\mathcal{T})}$ be a configuration drawn from the Ising model on $\mathcal{T}$ with external field $h$ and inverse temperature $s$ (and let $\theta=\tanh(s)$). Denote by $\E_h$ the expectation under the Ising Model with external field $h$. Then, if $\rho$ is the root of $\mathcal{T}$, $v_1,...,v_k$ are vertices of $\mathcal{T}$ and $\E_h(\sigma_\rho)=0$, the following inequality holds:
$$\E_h\left(\sigma_\rho|\sigma_{v_1}=\cdots=\sigma_{v_k}=+\right)\le 50\cdot\SL{i=1}^k\E_0\left(\sigma_\rho|\sigma_{v_i}=+\right)=50\cdot\SL{i=1}^k\theta^{\text{dist}(\rho,v_i)}.$$
\end{prop}
Before we prove Proposition \ref{main-prop}, we mention the following Theorem, which is Corollary 1.3 in \cite{D-S-S}.
\begin{thm}\label{D-S-S}
Let $g:V\to[-\infty,\infty]$ be an external field. Then, for all $u,v\in V$
$$\E_g(\sigma_u\sigma_v)-\E_g(\sigma_u)\E_g(\sigma_v)\le\E_0(\sigma_u\sigma_v).$$
\end{thm}
\begin{proof}[Proof of Proposition \ref{main-prop}]
In what follows, when we omit the external field in the expectation, we mean external field equal to $h$. Of course, we can write
\begin{align*}
\E\left(\sigma_\rho|\sigma_{v_1}=\cdots=\sigma_{v_k}=+\right)&=\SL{i=1}^k\left[\E\left(\sigma_\rho|\sigma_{v_1}=\cdots=\sigma_{v_i}=+\right)-\E\left(\sigma_\rho|\sigma_{v_1}=\cdots=\sigma_{v_{i-1}}=+\right)\right]\\&=\SL{i=1}^k\left[\E_{h_{i-1}}\left(\sigma_\rho|\sigma_{v_i}=+\right)-\E_{h_{i-1}}\left(\sigma_\rho\right)\right],
\end{align*}
where $h_{i-1}$ is the new external field coming from conditioning on $\{\sigma_{v_1}=\cdots=\sigma_{v_{i-1}}=+\}$ as well. We now calculate and bound each term of the above sum separately. For convenience, let $p_i=\PR_{h_{i-1}}(\sigma_{v_i}=+), m_i^\pm=\E_{h_{i-1}}(\sigma_\rho|\sigma_{v_i}=\pm)$ and observe that $m_i^+\ge m_i^-$. It is not hard to see that
$$p_i\ge\PR_0(\sigma_\rho=+|\sigma_{u_1}=\cdots=\sigma_{u_m}=-)\ge\dfrac{1}{1+e^{2d\beta}}\ge\dfrac{1}{100},$$
where $u_1,...,u_m$ (with $m\le d$) are the neighboring vertices of $\rho$.\\
Since $\E_{h_{i-1}}(\sigma_\rho)=p_im_i^++(1-p_i)m_i^-$, we can see that
\begin{align*}
\E_{h_{i-1}}(\sigma_\rho|\sigma_{v_i}=+)-\E_{h_{i-1}}(\sigma_\rho)=(1-p_i)(m_i^+-m_i^-)
\end{align*}
and
\begin{align*}
\text{Cov}(\sigma_\rho,\sigma_{v_i})&=\E(\sigma_\rho\sigma_{v_i})-\E(\sigma_\rho)\E(\sigma_{v_i})\\&=p_im_i^+-(1-p_i)m_i^--(2p_i-1)(p_im_i^++(1-p_i)m_i^-)\\&=2p_i(1-p_i)(m_i^+-m_i^-).
\end{align*}
It is now clear that
\begin{align*}
\E_{h_{i-1}}\left(\sigma_\rho|\sigma_{v_i}=+\right)-\E_{h_{i-1}}\left(\sigma_\rho\right)&\le50\cdot\text{Cov}_{h_{i-1}}(\sigma_\rho,\sigma_{v_i})\\&\le50\cdot\text{Cov}_0(\sigma_\rho,\sigma_{v_i})=50\cdot\E_0(\sigma_\rho|\sigma_{v_i}=+),
\end{align*}
where the last inequality comes from Theorem \ref{D-S-S}. The proof of the Proposition is complete.
\end{proof}
We wish to apply Proposition \ref{main-prop} to the tree of Self-Avoiding Walks of the graph $G$, which we define below (for more details see also \cite{Weitz}).
\begin{defn}
Let $G$ be a graph. On the vertices of $G$ fix an enumeration and let $v$ be one of the vertices of $G$. The tree of Self-Avoiding Walks $T_{\text{SAW}}(G,v)$ of $G$ rooted at $v$ is the tree of walks originating at $v$ that do not intersect themselves, except when the walk closes a cycle, at which point the walk ends. When that happens, the second copy of the twice-appearing vertex is fixed to be a $+$, if the edge closing the cycle is larger than the edge starting the cycle, and $-$ otherwise. In this case, comparing the edges really means comparing the respective vertices, with respect to the enumeration we initially fixed. It is easy to see that the degree of each vertex of $T_{\text{SAW}}(G,v)$ is at most $d$.
\end{defn}
In the above construction, there is an implied map $\phi$, which maps every vertex of $T_{\text{SAW}}(G,v)$ to its label in the original graph $G$. The reason why this construction is so useful, is the following Theorem, essentially proven in \cite{Weitz}.
\begin{thm}\label{Weitz}
Let $y\in V$, and consider the Ising model on $G$, conditioned on $\sigma_y=+$. Also, let $\sigma_{c,y}$ be the corresponding conditions on $T_{\text{SAW}}(G,v)$ obtained as follows: Every vertex $u$ on $T_{\text{SAW}}(G,v)$ which has $\phi(u)=y$ and has not been set to have a specific value from the construction of $T_{\text{SAW}}(G,v)$, is set to be a $+$. Then,
$$\PR_G(\sigma_v=+|\sigma_y=+)=\PR_{T_{\text{SAW}}(G,v)}(\sigma_\rho=+|\sigma_{c,y}).$$
\end{thm}
\begin{proof}[Proof of Lemma \ref{main-lem}]
Let $\mathcal{T}$ be the tree of the self avoiding paths on $G$ rooted at $x$ and $\phi$ the function that indicates the label of each vertex in $\mathcal{T}$. As in Theorem \ref{Weitz}, for a vertex $y$, set $\sigma_y=+$, to obtain conditions $\sigma_{c,y}$. Then,
$$\E(\sigma_x\sigma_y)=\E(\sigma_x|\sigma_y=+)=\E_{\mathcal{T}}(\sigma_\rho|\sigma_{c,y})=\E_{h,\mathcal{T}}(\sigma_\rho|\sigma_u=+,\ \forall u\in\phi^{-1}(y)),$$
where the external field $h$ comes from the conditions imposed by the construction of $\mathcal{T}$. Due to the symmetry of the Ising model, it is easy to see that $\E_h(\sigma)=0$. So, due to Proposition \ref{main-prop}, for any $y$
$$\E(\sigma_x\sigma_y)\le50\cdot\SL{u:\phi(u)=y}\theta^{\text{dist}(\rho,u)}.$$
Adding up all of those and using the fact that $\mathcal{T}$ is a subgraph of the infinite $d$-regular tree, we get that for every $x$,
\begin{align*}
\SL{y\in V(G)}\E(\sigma_x\sigma_y)\le1+50\cdot\SL{k=1}^\infty d(d-1)^{k-1}\theta^k\le\dfrac{50d}{d-1}\cdot\dfrac{1}{1-(d-1)\theta}\le\dfrac{100d}{(d-1)^2\cdot(\beta_c-s)}.
\end{align*}
We have proven the desired result.
\end{proof}

\end{document}